\documentclass[notitlepage,11pt]{extarticle}
\usepackage[T2A]{fontenc}
\usepackage{amssymb, amsmath}

\usepackage{geometry}
\usepackage[utf8x]{inputenc}
\usepackage{ucs}
\usepackage[english]{babel}
\usepackage[unicode,colorlinks=true,linkcolor=.]{hyperref}
\usepackage{pifont}
\usepackage{yfonts}
\usepackage{textcomp}
\usepackage{float}
\usepackage{mathtools}
\usepackage{amsmath}
\usepackage{amsfonts}
\usepackage{amsthm}
\usepackage{amssymb}
\usepackage{bm}
\usepackage{proof}
\usepackage{longtable}
\usepackage{pdfpages}

\newcommand{\qs}[1]{{\emph{``#1''}}}

\newcommand{\eb}[1]{{\left(#1\right)}}
\newcommand{\ecb}[1]{\left\{#1\right\}}
\newcommand{\esb}[1]{\left[#1\right]}

\newcommand{\CL}{\texttt{CL}}
\newcommand{\IL}{\texttt{IL}}
\newcommand{\BL}{\texttt{BL}}
\newcommand{\BSL}{\texttt{BSL}}
\newcommand{\ML}{\texttt{ML}}
\newcommand{\PEL}{\texttt{PEL}}
\newcommand{\PL}{\texttt{PL}}

\newcommand{\withor}{\esb{\lor}}

\newcommand{\comment}[1]{}

\newcommand{\equ}{\sim}
\newcommand{\xtophi}{{\esb{x_0/\phi}}}
\newcommand{\xtopsi}{{\esb{x_0/\psi}}}

\theoremstyle{definition}
\newtheorem{tDefinition}{Definition}[section]

\theoremstyle{plain}
\newtheorem{tTheorem}{Theorem}[section]
\newtheorem*{tMainTheorem}{Main theorem}
\newtheorem{tLemma}{Lemma}[section]
\newtheorem{tCorollary}{Corollary}[section]

\newcommand{\theorem}[3]{
	\begin{tTheorem}\label{#1}
		{#2}
	\end{tTheorem}
	\begin{proof}
		{#3}
	\end{proof}
}

\newcommand{\lemma}[3]{
	\begin{tLemma}\label{#1}
		{#2}
	\end{tLemma}
	\begin{proof}
		{#3}
	\end{proof}
}

\newcommand{\corollary}[3]{
	\begin{tCorollary}\label{#1}
		{#2}
	\end{tCorollary}
	\begin{proof}
		{#3}
	\end{proof}
}

\newcommand{\definition}[2]{
	\begin{tDefinition}\label{#1}
		{#2}
	\end{tDefinition}
}

\usepackage{tocloft}

\begin{document}

\begingroup
	\let\clearpage\relax
	\title{On Computationally Efficient Subsystems of Propositional Logic}
	\date{February 7, 2016}
	\author{
		Inga Lev
	}
	
	\maketitle
	\begin{abstract}
	
	In this paper, we show that the derivability problem for the primal propositional logic
	remains solvable in polynomial time upon adding a certain form of the principle of equivalent form
	substitution; and that, upon adding another form of this principle, it becomes co-NP-hard.
	
	\end{abstract}
	
	\setlength{\cftbeforesecskip}{0.2em}
	\tableofcontents
\endgroup

\newpage

\section{Introduction}

There are numerous real-world tasks that include determining whether a given formula is derivable
from given assumptions in a certain logic. One possible area is access control; the task of determining
whether certain access is to be granted based on some prior knowledge could be reduced to the task of
determining whether the proposition “Access should be granted” is derivable from the corresponding
assumptions.

Obviously, for such a system to be usable, the corresponding logic should conform to certain
desirable requirements, and the derivability problem should be reasonably easy (when expressed in
terms of the sequent length).

It is known that the derivability problem for classical propositional logic is co-NP-complete; and,
by a result of Statman \cite{RS1979}, intuitionistic logic is PSPACE-complete. Additionally, it is easy to
prove that any intermediate logic is co-NP-hard.

Y. Gurevich and I. Neeman proved, in the context of their research into infon logic and Distributed
Knowledge Authorization Language, that the derivability problem is solvable in polynomial time for
a certain subsystem of disjunction-free minimal logic, called \emph{primal logic} ($\PL$), which is, in a certain
sense, the smallest useful logic.

However, $\PL$ is so limited that it lacks certain properties one expects from the real-world logic. In
particular, it does not conform to the \emph{principle of equivalent formula substitution}; that is, it is possible
to turn a derivable formula into a non-derivable one or vice versa by changing one of its subformulas
for an equivalent one (we consider two formulas $\phi$ and $\psi$ to be equivalent if they are derivable from
each other, $\phi \vdash \psi$, and $\psi \vdash \phi$).

The main objective of this paper is to show that, for every logic that extends $\PL$ and conforms to
the strong form of principle of equivalent formula substitution, the derivability problem is co-NP-hard;
yet that there is a logic that extends $\PL$ and conforms to the week form of principle of equivalent
formula substitution, for which the derivability problem is solvable in cubic time. It is achieved by
showing that $\PL$ with the strong form of principle of equivalent formula substitution could be used to
solve the derivability problem for classical logic, while $\PL$ itself could be used to solve the derivability
problem for $\PL$ with the weak form of this principle.

Additionally, we will show that the smallest logic containing the strong form of the principle of
equivalent substitution is different from minimal logic.

This work was partially done during my study at the Faculty of Mathematics, Higher School of Economics, Moscow, Russia.
I am grateful to my supervisor Lev Beklemishev for suggesting me to work on this problem and for his guidance and advice.

\section{Preliminaries}
In this paper we study propositional logics in sequential format. We will mostly concentrate on the
disjunction-free language.

Let us define some basic notions used in this paper.

\subsection{Logics}

We define basic logic, denoted as $\BL$, by the following rules:

\begin{longtable}{c c}
$\top$ & ${\vdash \top}$ \\\\

$\texttt{x2x}$ & $\phi \vdash \phi$ \\\\

$\texttt{Premise inflation}$ & $\infer{\Gamma, \phi \vdash \psi}{\Gamma \vdash \psi}$ \\\\

$\texttt{Cut}$ & $\infer{\Gamma \vdash \psi}{\Gamma \vdash \phi &\quad \Gamma, \phi \vdash \psi}$ \\\\

$\land \texttt{E}_l, \land \texttt{E}_r$ &\quad $\infer{\Gamma \vdash \phi}{\Gamma \vdash \eb{\phi \land \psi}},\quad \infer{\Gamma \vdash \psi}{\Gamma \vdash \eb{\phi \land \psi}}$ \\\\

$\land \texttt{I}$ &\quad $\infer{\Gamma \vdash \eb{\phi \land \psi}}{\Gamma \vdash \phi &\quad \Gamma \vdash \psi}$ \\\\
\end{longtable}

We define primal logic, denoted as $\PL$, by adding the following inference rules to $\BL$:

\begin{longtable}{c c}
$\to \texttt{E}$ & $\infer{\Gamma \vdash \psi}{\Gamma \vdash \phi &\quad \Gamma \vdash \eb{\phi \to \psi}}$ \\\\

$\to \texttt{IW}$ & $\infer{\Gamma \vdash \eb{\phi \to \psi}}{\Gamma \vdash \psi}$
\end{longtable}

In \cite{G205}, Gurevich and Beklemishev prove that the multiple derivability problem for $\PL$ is solvable
in linear-time, while the same problem for $\PL\withor$ is co-NP-hard.

We define minimal logic, denoted as $\ML$, by adding the following inference rule to $\PL$:

\begin{longtable}{c c}
$\to\texttt{I}$ & $\infer{\Gamma \vdash \eb{\phi \to \psi}}{\Gamma, \phi \vdash \psi}$
\end{longtable}

Note that, by premise inflation, it extends inference rule $\to \texttt{IW}$.

We define intuitionistic logic, denoted as $\IL$, is obtained from $\ML$ by adding the $\bot$ axiom: $\bot \vdash \phi$.

It was shown in \cite{RS1979} that the derivability problem for $\IL\withor$ is PSPACE-complete.

We define classical logic, denoted as $\CL$, by adding the following inference rule to $\IL$:

\begin{longtable}{c c}
$\genfrac{}{}{0pt}{0}{\texttt{Disjunction-free law}}{\texttt{of excluded middle}}$ & $\infer{\Gamma \vdash \psi}{\Gamma, \phi \vdash \psi &\quad \Gamma, \eb{\phi \to \bot} \vdash \psi}$
\end{longtable}

It is well known that the derivability problem for $\CL\withor$ is co-NP-complete.

\subsection{Equivalent sequents}

\definition{equivalence}
{
	We say that $\phi$ is \emph{equivalent} to $\psi$, or that $\phi \equ \psi$ in a given logic, if both $\phi \vdash \psi$ and $\psi \vdash \phi$ are theorems in that logic.

	Additionally, we say that $\phi$ is equivalent to $\psi$ under assumptions $\Gamma$, or that $\phi \equ_\Gamma \psi$, if both $\Gamma, \phi \vdash \psi$ and $\Gamma, \psi \vdash \phi$ are theorems.
	It is a weaker relation, as for every $\phi$ and $\psi$, the relation $\phi \equ \psi$ implies $\phi \equ_\Gamma \psi$ for every $\Gamma$.

	Note that both forms of equivalence relation are reflexive, symmetric, and transitive. That is, $\phi \equ_\Gamma \phi$ are theorems for every $\phi$; and, if $\phi \equ_\Gamma \psi$ and $\psi \equ_\Gamma \chi$ are all theorems for some $\phi$, $\psi$, and $\chi$, then $\psi \equ_\Gamma \phi$ and $\phi \equ_\Gamma \chi$ are theorems too. The same applies to the simple $\equ$ relation.
}

\subsection{Polynomial-time locality}

A logic is said to have the \emph{polynomial-time locality} property, if there is a polynomial-time algorithm
which, given the theorem $\sigma$, constructs a set of formulas $\phi_i$ such that there exists a proof of $\sigma$ in the
logic concerned such that every step of the proof contains only the $\phi_i$ formulas.

Related to this concept is the \emph{subformula} property. A logic is said to have the subformula property
if for any theorem $\sigma$ there exists a proof such that every sequent of it contains only the subformulas
of $\sigma$.

The subformula property implies the polynomial-time locality property, as any sequent $\sigma$ of the
length $n$ has at most $n$ subformulas, which all could be obtained in a linear time.

As noted in \cite{G198} and \cite{G205}, $\ML$ and $\ML\withor$ have the subformula property.

\subsection{Ordering of logics}

\definition{derivedRule}
{
	The inference rule $R$ is said to be a derived rule in logic $A$,  if it could be described
as a combination of some inference rules of $A$.
}

We say that $A \le B$ for logics $A$ and $B$, if every inference rule of $A$ is a derived rule in $B$.

Note that this implies that all theorems in $A$ are theorems in $B$ as well, as their proofs in $A$
could be translated to $B$ (in linear time with respect to the length of the proof). Set of all theorems
in $A$ is the subset of the set of all theorems in $B$.

Obviously, if logic $B$ is obtained from logic $A$ by adding some inference rule to it, then $A \le B$.

We say that logics $A$ and $B$ are equivalent ($A \equiv B$), if $A \leq B$ and $B \leq A$. Set of all theorems in $A$
is then equal to the set of all theorems in $B$, so $A$ and $B$ are of the same complexity, and the proof of any
theorem in $A$ could be translated to $B$ in linear time and vice versa

We say that $A < B$, if $A \le B$ and $A \not\equiv B$. Note that set of all theorems in $A$ may still be equal to
the set of all theorems in $B$.

\subsection{Complexity of logic}

We define the complexity of logic as the complexity of the algorithmic problem to decide whether a given sequent is a theorem in the logic concerned or not.

For example, we say that the complexity of the logic $A$ is polynomial-time, if there are number $c$ and $k$ and an algorithm $F$ such that, given a sequent $\sigma$ of the length $n$, $F$ terminates after no more than $c*n^k$ steps and correctly decides whether $\sigma$ is a theorem in $A$ or not.

\subsection{Compatibility}

We say that, in context of the logic $B$, the sequent $\sigma$ is $A$-compatible, iff it is simultaneously a theorem,
or a non-theorem, in both $B$ and $A$.

If $B > A$, then the set of all $A$-compatible sequents in logic $B$ is equal to the disjoint union of the
sets of all theorems in $A$ and all non-theorems in $B$.

\subsection{Reduction}

We say that logic $B$ is reducible to logic $A$, if there is a theoremhood-preserving mapping (reduction)
from sequents in $B$ to $A$-compatible sequents.

If $t\eb{\sigma}$ is such a mapping, $p_t\eb{\sigma}$ is the time required to compute $t\eb{\sigma}$, and $p_A\eb{\sigma}$ is the time required
to compute the derivability problem for $\sigma$ in $A$, then, obviously, $p_B\eb{\sigma} \leq p_t\eb{\sigma}p_A\eb{t\eb{\sigma}}$.

Now, if $t'\eb{n}$ is the maximum length of all $t\eb{\sigma}$ for all $\sigma$ with the length no more than $n$; $p'_t\eb{n}$ is the maximum time required to compute $t\eb{\sigma}$ for any such $\sigma$; and $p'_A\eb{n}$ is the maximum time
required to compute the derivability problem in $A$ for any such $\sigma$ (that is, the complexity of $A$), then
$p'_B\eb{n} \leq p'_t\eb{n}p'_A\eb{t'\eb{n}}$. If $t$ does not increase sequent length too much, and does not require too much
time to compute, then $B$ could not be significantly more complex than $A$. In particular, if $t$ is in P,
and $B$ is NP-hard, then $A$ is NP-hard too.

In \cite{G205}, Gurevich and Beklemishev prove that $\CL\withor$ is reducible to $\PL\withor$, showing that it
is co-NP-hard; in \cite{YS2009}, Savateev proves that $\PL$ is reducible to $\BL$, againg showing that $\PL$ is
polynomial-time decidable.

\lemma{middleLemma}
{
If $A \le B \le C$ and $C$ is reducible to $A$ by a mapping $t$, then $C$ is reducible to $B$ by the same mapping.
}
{
Every $A$-compatible sequent in $C$ is $B$-compatible as well.
}

\section{The principle of equivalent formula substitution}

The \emph{principle of compositionality} states that \qs{the meaning of a complex expression is fully determined
by its structure and the meanings of its constituents} (\cite{SEP}).

This property is, understandably, desirable in the logics used in the real world.

We will study a similar property in this paper:

\definition{equivalentFormulaSubstitution}
{
The logic $A$ is said to satisfy the \emph{principle of equivalent formula substitution}, if every of the following family of inference rules, generated by all possible expressions $F$, is a derived rule in $A$:

\begin{longtable}{c c}
$\texttt{E}_F$ & $\infer{\Gamma, F\xtophi \vdash F\xtopsi}{\Gamma, \phi \vdash \psi &\quad \Gamma, \psi \vdash \phi}$
\end{longtable}

That is, if two formulas are equivalent under assumptions $\Gamma$, then they yield equivalent formulas
under assumptions $\Gamma$ when substituted in any expression.

The weak form of this principle is defined by using the following family of inference rules instead:

\begin{longtable}{c c}
$\texttt{E}_{0,F}$ & $\infer{F\xtophi \vdash F\xtopsi}{\phi \vdash \psi &\quad \psi \vdash \phi}$,
\end{longtable}
}

That is, if two formulas are equivalent, then they yield equivalent formulas when substituted in
any expression.

Let us define the following inference rules by substituting $\eb{x_0 \to \chi}$ and $\eb{\chi \to x_0}$ in place of $F$:

\begin{longtable}{c c}
$\texttt{E1}$ & $\infer{\Gamma, \eb{\phi \to \chi} \vdash \eb{\psi \to \chi}}{\Gamma, \phi \vdash \psi &\quad \Gamma, \psi \vdash \phi}$ \\\\
$\texttt{E2}$ & $\infer{\Gamma, \eb{\chi \to \phi} \vdash \eb{\chi \to \psi}}{\Gamma, \phi \vdash \psi &\quad \Gamma, \psi \vdash \phi}$ \\\\
$\texttt{E1}_0$ & $\infer{\eb{\phi \to \chi} \vdash \eb{\psi \to \chi}}{\phi \vdash \psi &\quad \psi \vdash \phi}$ \\\\
$\texttt{E2}_0$ & $\infer{\eb{\chi \to \phi} \vdash \eb{\chi \to \psi}}{\phi \vdash \psi &\quad \psi \vdash \phi}$
\end{longtable}

By extending $\PL$ with the rules $\texttt{E1}$, $\texttt{E2}$, or both, we obtain logics we will refer to as $\PEL1$, $\PEL2$, and
$\PEL$, respectively. We similarly define smaller logics $\PEL1_0$, $\PEL2_0$, and $\PEL_0$.

\lemma{pel0clintermediate}
{
	If $\PEL_0 \leq A \leq \CL$, then $A$ satisfies the weak form of the principle of equivalent formula
	substitution.

	If $\PEL \leq A \leq \CL$, then $A$ satisfies the principle of equivalent formula substitution.
}
{
	We will prove the claim for the full form of the principle of equivalent formula substitution;
	proof for the weak form is analogous.

	Let us prove this claim by the induction on the length of $F$.

	First of all, note that the claim is obvious if $F$ is invariant under replacement of $x_0$ to $\phi$ and $\psi$
	(that is, $F$ does not contain $x_0$ as its part), since $E_F$ is then actually $\infer{\Gamma, \chi \vdash \chi}{\Gamma, \phi \vdash \psi &\quad \Gamma, \psi \vdash \phi}$, which is a
	derived rule in $\BSL$ from $\texttt{x2x}$ and $\texttt{Premise inflation}$.

	The claim is true if $F$ is a placeholder or a constant: If it is $x_0$, then $E_F$ is actually $\infer{\Gamma, \phi \vdash \psi}{\Gamma, \phi \vdash \psi &\quad \Gamma, \psi \vdash \phi}$,
	which is a tautology. And if it is not $x_0$, $E_F$ is also a tautology.

	Now, let us suppose that $F$ is a combination of two shorter expressions (say, $G$ and $H$) and one binary operator. By the induction hypothesis we have that both
	$\infer{\Gamma, G\xtophi \vdash G\xtopsi}{\Gamma, \phi \vdash \psi &\quad \Gamma, \psi \vdash \phi}$ and
	$\infer{\Gamma, H\xtophi \vdash H\xtopsi}{\Gamma, \phi \vdash \psi &\quad \Gamma, \psi \vdash \phi}$ are derived rules; so, assuming that we have $\phi \equ_\Gamma \psi$ in the derivation at this moment, we get both
	$\Gamma, G\xtophi \vdash G\xtopsi$ and $\Gamma, H\xtophi \vdash H\xtopsi$ by the induction hypothesis. There are two possibilities:

	\begin{itemize}
	\item{
		$F = \eb{G \land H}$.
	
		By the $\texttt{x2x}$, $\land\texttt{E}$, and $\texttt{Premise inflation}$ rules we obtain $\Gamma, F\xtophi \vdash G\xtophi$ and $\Gamma, F\xtophi \vdash H\xtophi$.
	
		Then, with the $\texttt{Premise inflation}$ and $\texttt{Cut}$ rules we obtain $\Gamma, F\xtophi \vdash G\xtopsi$ and $\Gamma, F\xtophi \vdash H\xtopsi$.
	
		Finally, by the $\land\texttt{I}$ rule we obtain $\Gamma, F\xtophi \vdash F\xtopsi$.
	}
	
	\item{
		$F = \eb{G \to H}$.
	
		By the $\texttt{E}1$ rule we obtain $\Gamma, \eb{G\xtophi \to H\xtophi} \vdash \eb{G\xtopsi \to H\xtophi}$.
	
		By the $\texttt{E}2$ rule we obtain $\Gamma, \eb{G\xtopsi \to H\xtophi} \vdash \eb{G\xtopsi \to H\xtopsi}$.
	
		Combining these two by the $\texttt{Cut}$ rule, we obtain $\Gamma, \eb{G\xtophi \to H\xtophi} \vdash \eb{G\xtopsi \to H\xtopsi}$.
	}
	\end{itemize}
	
	Thus, $\Gamma, F\xtophi \vdash F\xtopsi$.
}

Now it is clear that $\PEL$ is the minimal logic which is larger than $\PL$ and satisfies the principle of
equivalent formula substitution; and, similarly, $\PEL_0$ is the minimal logic which is larger than $\PL$ and
satisfies the weak form of the principle of equivalent formula substitution.

\section{Models}

\subsection{The degenerate implication}

Let us define the following inference rule of degenerate implication:

\begin{longtable}{c c}
$\to\texttt{ED}$ & $\infer{\Gamma \vdash \psi}{\Gamma \vdash \eb{\phi \to \psi}}$
\end{longtable}

In logics with such a rule, $\eb{\phi \to \psi} \equ \psi$.

Note that $\to\texttt{ED}$ extends inference rule $\to\texttt{E}$, in the sense that it allows us to infer everything $\to\texttt{E}$
does. Additionally, note that it extends all of the $\texttt{E}1$, $\texttt{E}2$, $\texttt{E}1_0$, $\texttt{E}2_0$ rules.

The \emph{valuation} $V$ is any function from the set of all sequents into $\ecb{\top, \bot}$ satisfying the following
conditions:

\begin{itemize}
\item{$V\eb{\vdash \top}=\top$.}
\item{$V\eb{\vdash \eb{\phi \land \psi}}=\top$ iff $V\eb{\vdash \phi}=V\eb{\vdash \psi}=\top$.}
\item{$V\eb{\vdash \eb{\phi \to \psi}}=V\eb{\vdash \psi}$}
\item{$V\eb{\Gamma \vdash \phi}=\top$ is true iff $V\eb{\vdash \Gamma_i}=\bot$ for some $\Gamma_i \in \Gamma$, or $V\eb{\vdash \phi}=\top$.}
\end{itemize}

It is easy to see that the valuation is completely defined by the images of variables.

We won’t prove the completeness theorem for such a model; instead, we will prove the soundness
theorem.

\lemma{pdilCorrectness}
{
	If $\sigma$ is a theorem in $\PL + \to\texttt{ED}$, then $V\eb{\sigma}=\top$ for every valuation $V$.
}
{
Let us prove this by buildup of the set of all theorems in such a logic. That is, we need to
check, that, if some sequent could be obtained from other sequents by applying some inference rule,
and each of these other sequents is already evaluated as $\top$, then the resulting sequent should also be
evaluated as $\top$.

\begin{description}
\item[$\top$:]{
	$V\eb{\vdash \top}=\top$.
}
\item[$\texttt{x2x}$:]{
	For every $\phi$, $V\eb{\phi \vdash \phi}=\top$ independent of $V\eb{\phi}$ value.
}
\item[$\texttt{Premise inflation}$:]{
	For every $\Gamma$, $\phi$, and $\psi$, if $V\eb{\Gamma \vdash \phi}=\top$, then $V\eb{\Gamma, \psi \vdash \phi}=\top$.
}
\item[$\texttt{Cut}$:]{
	For every $\Gamma$, $\phi$, and $\psi$, if $V\eb{\Gamma \vdash \phi}=\top$, then either (1) $V\eb{\vdash \Gamma_i} = \bot$ for some $\Gamma_i$ in $\Gamma$, or (2) $V\eb{\vdash \phi} = V\eb{\vdash \Gamma_i} = \top$ for every $\Gamma_i$ in $\Gamma$.

	If the first case, for every $\chi$, $V\eb{\Gamma \vdash \chi} = \top$; so, in particular, $V\eb{\Gamma \vdash \psi} = \top$ (for the sake of
simplicity, we will omit this case for the similar inference rules considered below, and will only
consider the case of $V\eb{\vdash \Gamma_i}$ being equal to $\top$ for every $\Gamma_i$
in $\Gamma$).

	In the second case, $V\eb{\Gamma, \phi \vdash \psi} = \top$ tells us that $V\eb{\Gamma \vdash \psi} = V\eb{\vdash \psi} = \top$
}
\item[$\land\texttt{E}_l, \land\texttt{E}_r$:]{
	For every $\Gamma$, $\phi$, and $\psi$, if $V\eb{\Gamma \vdash \eb{\phi \land \psi}} = \top$, then, omitting the case of some $\Gamma_i$ being
	evaluated to $\bot$, we obtain that $V\eb{\vdash \eb{\phi \land \psi}} = \top$. So, $V\eb{\vdash \phi} = V\eb{\vdash \psi} = \top$; $V\eb{\Gamma \vdash \phi} = V\eb{\Gamma \vdash \psi} = \top$.
}
\item[$\land\texttt{I}$:]{
	For every $\Gamma$, $\phi$, and $\psi$, if $V\eb{\Gamma \vdash \phi} = V\eb{\Gamma \vdash \psi} = \top$, then, omitting the case of some $\Gamma_i$ being
	evaluated to $\bot$, we obtain that $V\eb{\vdash \phi} = V\eb{\vdash \psi} = \top$. So, $V\eb{\Gamma \vdash \eb{\phi \land \psi}} = V\eb{\vdash \eb{\phi \land \psi}} = \top$.
}
\item[$\to\texttt{ED}$:]{
	For every $\Gamma$, $\phi$, and $\psi$, if $V\eb{\Gamma \vdash \eb{\phi \to \psi}} = \top$ then, omitting the case of some $\Gamma_i$ being
	evaluated to $\bot$, we obtain that $V\eb{\Gamma \vdash \psi} = V\eb{\vdash \psi} = \top$.
}
\item[$\to\texttt{E}$:]{
	Is a particular case of $\to\texttt{ED}$.
}
\item[$\to\texttt{IW}$:]{
	For every $\Gamma$, $\phi$, and $\psi$, if $V\eb{\Gamma \vdash \psi} = \top$ then, omitting the case of some $\Gamma_i$ being evaluated to
	$\bot$, we obtain that $V\eb{\Gamma \vdash \eb{\phi \to \psi}} = V\eb{\vdash \eb{\phi \to \psi}} = V\eb{\vdash \psi} = \top$.
}
\end{description}
}

\subsection{The models of PL}

Let us define the Kripke models in the way similar to the one for intuitionistic logic, with the only
exception of how we treat the $\to$ relation. In intuitionistic logic, $\phi \to \psi$ is evaluated as true in the
world $a$ iff in every world $b \geq a$ where $\phi$ is evaluated as true, $\psi$ is evaluated as true as well. We instead
define $x \vdash y$ in such a fashion (or, more formally, $\Gamma \vdash y$ is evaluated as true in world $a$ iff for every
world $b \geq a$ where all the formulas of $\Gamma$ are evaluated as true, $y$ is evaluated as true as well), and we
allow for $x \to y$ to be evaluated in any way as long as it satisfies two conditions:

\begin{itemize}
\item{Wherever $y$ is evaluated as true, $x \to y$ also has to be evaluated as true.}
\item{Wherever $x$ and $x \to y$ are both evaluated as true, $y$ also has to be evaluated as true.}
\end{itemize}

Both completeness and soundness theorems for such a model were proved in \cite{G198}.

\section{Certain observations on logics}

\lemma{pdilIsNotMlOrCl}
{
	$\ML \not< \PL+\to\texttt{ED}$, $\PL+\to\texttt{ED} \not< \CL$. That is, there are sequents which are theorems in $\ML$
but are not theorems in $\PL+\to\texttt{ED}$; and there are sequents which are theorems in $\PL+\to\texttt{ED}$ but are
not theorems in $\CL$.
}
{
	The sequent $\eb{x_1 \to x_2} \vdash x_2$ is not a theorem in $\CL$ (as one could see by evaluating both $x_1$ and
	$x_2$ as $\bot$). However, it is a theorem in $\PL+\to\texttt{ED}$, obtained from $\eb{x_1 \to x_2} \vdash \eb{x_1 \to x_2}$ by $\to\texttt{ED}$ rule.

	Let us define the valuation $V$ by letting $V\eb{x_i} = \bot$ for every $i$. Obviously, $V\eb{\vdash x_1}$ is false. By Lemma \ref{pdilCorrectness}, $\vdash x_1$ is not a theorem in $\PL+\to\texttt{ED}$.
}

\lemma{pelIsNotMl}
{
	$\PEL \ne \ML$.
}
{
	Immediately follows from lemma \ref{pdilIsNotMlOrCl}.
}

\lemma{pelIsNotPl}
{
	$\PEL1_0, \PEL2_0 \not= \PL$.
}
{
	First, note that $\eb{x \to x} \vdash \eb{x \to x}$, $x \vdash \eb{x \land x}$, and $\eb{x \land x} \vdash x$ are all theorems in $\PL$

	Now, let us define the Kripke model consisting of two worlds embedded in each other such that:

	\begin{itemize}
	\item{$\vdash x_1$ (and thus $\vdash \eb{x_1 \land x_1}$ as well) is only evaluated as true in the upper world}
	\item{$\vdash \eb{x_1 \to x_1}$ is evaluated as true in both worlds}
	\item{$\vdash \eb{x_1 \to \eb{x_1 \land x_1}}$ and $\eb{\eb{x \land x} \to x}$ are only evaluated as true in the upper world}
	\end{itemize}

	Such valuation of the $\to$ relation satisfies the definition of Kripke model.

	In this model, neither $\eb{x \to x} \vdash \eb{\eb{x \land x} \to x}$ nor $\eb{x \to x} \vdash \eb{x \to \eb{x \land x}}$ are evaluated as true
	in the lower world. Thus, neither is a theorem in $\PL$.

	Yet, $\eb{x \to x} \vdash \eb{\eb{x \land x} \to x}$ is a theorem in $\PEL1_0$; and $\eb{x \to x} \vdash \eb{x \to \eb{x \land x}}$ is a
	theorem in $\PEL2_0$ (and $\PEL2_0$).
}

\section{The complexity of PEL}

\subsection{The complexity of CL}

\lemma{clComplexity}
{
	$\CL\esb{\lor}$ is reducible to $\CL$ by a polynomial-time mapping.
}
{
	Let us define the mapping $t$ on both formulas and sequents as follows:

	\begin{itemize}
		\item{$t\eb{x} = x$, where $x$ is a variable or constant}
		\item{$t\eb{\phi \land \psi} = \eb{t\eb{\phi} \land t\eb{\psi}}$}
		\item{$t\eb{\phi \lor \psi} = \eb{\eb{t\eb{\phi} \to \bot} \to \eb{\eb{t\eb{\psi} \to \bot} \to \bot}}$}
		\item{$t\eb{\phi \to \psi} = \eb{t\eb{\phi} \to t\eb{\psi}}$}
		\item{$t\eb{\phi_1, \ldots, \phi_n \vdash \psi} = t\eb{\phi_1}, \ldots, t\eb{\phi_n} \vdash t\eb{\psi}$}
	\end{itemize}

	Note that $t$ is idempotent.

	We will first prove that $\eb{\phi \lor \psi} \equ \eb{\eb{\phi \to \bot} \to \eb{\eb{\psi \to \bot} \to \bot}}$ in $\CL\esb{\lor}$. We get:

	$$\eb{\phi \lor \psi},\eb{\phi \to \bot},\eb{\psi \to \bot}, \phi \vdash \bot \quad \eb{\texttt{x2x}, \to \texttt{E}}$$
	$$\eb{\phi \lor \psi},\eb{\phi \to \bot},\eb{\psi \to \bot}, \psi \vdash \bot \quad \eb{\texttt{x2x}, \to \texttt{E}}$$
	$$\eb{\phi \lor \psi},\eb{\phi \to \bot},\eb{\psi \to \bot} \vdash \bot \quad \eb{\lor\texttt{I}}$$
	$$\eb{\phi \lor \psi},\eb{\phi \to \bot} \vdash \eb{\eb{\psi \to \bot} \to \bot} \quad \eb{\to I}$$
	$$\eb{\phi \lor \psi} \vdash \eb{\eb{\phi \to \bot} \to \eb{\eb{\psi \to \bot} \to \bot}} \quad \eb{\to I}$$

	On the other hand, we get:

	$$\eb{\eb{\phi \to \bot} \to \eb{\eb{\psi \to \bot} \to \bot}},\eb{\phi \to \bot},\eb{\psi \to \bot} \vdash \bot \quad \eb{\texttt{x2x}, \to \texttt{E}}$$
	$$\eb{\eb{\phi \to \bot} \to \eb{\eb{\psi \to \bot} \to \bot}},\eb{\phi \to \bot},\eb{\psi \to \bot} \vdash \eb{\phi \lor \psi} \quad \eb{\IL, \to \texttt{Cut}}$$
	$$\eb{\eb{\phi \to \bot} \to \eb{\eb{\psi \to \bot} \to \bot}},\eb{\phi \to \bot}, \psi \vdash \eb{\phi \lor \psi} \quad \eb{\texttt{x2x}, \lor \texttt{I}_r}$$
	$$\eb{\eb{\phi \to \bot} \to \eb{\eb{\psi \to \bot} \to \bot}},\eb{\phi \to \bot} \vdash \eb{\phi \lor \psi} \quad \eb{\texttt{Law of excluded middle}}$$
	$$\eb{\eb{\phi \to \bot} \to \eb{\eb{\psi \to \bot} \to \bot}}, \phi \vdash \eb{\phi \lor \psi} \quad \eb{\texttt{x2x}, \lor\texttt{I}_l}$$
	$$\eb{\eb{\phi \to \bot} \to \eb{\eb{\psi \to \bot} \to \bot}} \vdash \eb{\phi \lor \psi} \quad \eb{\texttt{Law of excluded middle}}$$

	By the induction on the formula length, we get that, in $\CL\esb{\lor}$, $\phi \equ t\eb{\phi}$.

	This means that $t$ is an equivalent formula substitution. From this immediately follows that $t$ is
	theoremhood-preserving.

	Now, let us suppose that in the image of mapping there are theorems that are not $\CL$-compatible.
	Let us take simplest of these theorems (in terms of the shortest derivation length in $\CL\esb{\lor}$), $\sigma$. All the
	theorems in that derivation are $\CL$-compatible under the mapping (otherwise, it would not be simplest).

	Let us consider the final step of that derivation. It could not involve any inference rule which does
	not deal with $\lor$ directly, as in that case, we could just apply the mapping $t$ to both the premises and
	conclusion, and obtain the derivation that is $\CL$-compatible. It could not be $\lor\texttt{I}_l$ or $\lor\texttt{I}_r$, since $\sigma$ is the
	result of the mapping, and does not contain $\lor$. Therefore, it should be $\lor\texttt{E}$.

	Then, for some $\Gamma, \phi, \psi, \chi$, the theorem $\sigma$ is equal to $\Gamma \vdash \chi$, and is obtained (in $\CL\esb{\lor}$) from theorems
	$\Gamma, \phi \vdash \chi$; $\Gamma, \psi \vdash \chi$; and $\Gamma \vdash \eb{\phi \lor \psi}$. By the choice of $\sigma$, both $\Gamma$ and $\chi$ do not contain disjunction;
	additionally, $t\eb{\Gamma \vdash \eb{\phi \lor \psi}} = \Gamma \vdash \eb{\eb{t\eb{\phi} \to \bot} \to \eb{\eb{t\eb{\psi} \to \bot} \to \bot}}$, $t\eb{\Gamma, \phi \vdash \chi} = \Gamma, t\eb{\phi} \vdash \chi$, and
	$t\eb{\Gamma, \psi \vdash \chi} = \Gamma, t\eb{\psi} \vdash \chi$ are all theorems in $\CL$. Yet these are enough to obtain $\sigma$ in $\CL$:
	
	$$\Gamma,\eb{t\eb{\phi} \to \bot},\eb{t\eb{\psi} \to \bot} \vdash \bot$$
	$$\Gamma,\eb{t\eb{\phi} \to \bot},\eb{t\eb{\psi} \to \bot} \vdash \chi$$
	$$\Gamma,\eb{t\eb{\phi} \to \bot}, t\eb{\psi} \vdash \chi$$
	$$\Gamma,\eb{t\eb{\phi} \to \bot} \vdash \chi$$
	$$\Gamma, \phi \vdash \chi$$
	$$\Gamma \vdash \chi$$

	Therefore, our assumption was wrong, and there is no such sequent $\sigma$; all the resulting sequents
	are $\CL$-compatible.
}

	Note that there is one additional consequence: every disjunction-free sequent of $\CL\esb{\lor}$ is $\CL$-compatible;
	disjunction-free theorems in $\CL$ do not require disjunction to prove them. Therefore,
	$\CL$ is reducible to $\CL\esb{\lor}$ with respect to a trivial linear-time mapping, which preserves all sequents that
	do not contain $\lor$ intact, and translates all other sequents to $\vdash \bot$. From this, we immediately obtain
	that $\CL$ is co-NP-complete.

\subsection{The complexity of ML}

\lemma{ilReducibleToMl}
{
	$\IL$ is reducible to $\ML$ with respect to a linear-time mapping.
}
{
	Let us show that if sequent $\sigma$ of $\IL$ contains helper formula ($\bot \to x_i$) as its antecedent for every
	$x_i$ mentioned in $\sigma$, then it is $\ML$-compatible; that is, it is either a theorem in $\ML$ or not a theorem in
	$\IL$.

	Let us assume that $\sigma$ is not a theorem in $\ML$. As it is not a theorem in $\ML$, there should be a world
	in some Kripke model such that, while all antecedents of $\sigma$ (including all $\eb{\bot \to x_i}$) are evaluated as
	true, its consequent is evaluated as false. Since the consequent of $\sigma$ is evaluated as false, some of the
	constants or variables contained in $\sigma$ should also be evaluated as false in this world. Therefore, as
	for every variable $x_i$ mentioned in $\sigma$, $\eb{\bot \to x_i}$ is evaluated as true, $\bot$ should be evaluated as false in
	this world. It is enough to use this world with everything above it as a counter-example in $\IL$ (with
	$\bot$ variable mapped to $\bot$ constant in $\IL$, for it is false in every world in this newly created model).
	Therefore, $\sigma$ is not a theorem in $\IL$.

	Now, let us consider a mapping that turns sequent $\sigma = \Gamma \vdash \phi$ to the sequent $\Gamma, \Delta \vdash \phi$, where $\Delta$ is
	the list of all helper formulas for $\sigma$. Such a mapping obviously turns theorems into theorems in every
	logic with $\texttt{Premise inflation}$ rule. It turns non-theorems into non-theorems in $\IL$, since it only adds
	as antecedents the consequent parts of antecedent-less theorems. As was shown above, it maps $\IL$
	sequents to $\ML$-compatible ones.

	Additionally, it takes linear time to obtain the resulting sequent.
}

	Now it immediately follows that $\ML$ is PSPACE-hard.

\subsection{The complexity of PEL}

Note that the following inference rules are derived rules in $\PL$:

$$\infer{\Gamma, \phi \vdash \eb{\phi \land \psi}}{\Gamma, \phi \vdash \psi},\quad \infer{\Gamma, \eb{\phi \land \psi} \vdash \phi}{\Gamma, \phi \vdash \psi}$$

\lemma{mlReducibleToPel1}
{
	$\ML$ is reducible to $\PEL1$ with respect to a polynomial-time mapping.
}
{
	Let us show that if sequent $\sigma$ of $\ML$ contains two \emph{helper formulas}

	$$\eb{\eb{\psi \land \omega} \to \eb{\psi \land \omega}}, \quad \eb{\eb{\psi \to \psi \land \omega} \to \eb{\psi \to \omega}}$$

	in its antecedent for every pair of \emph{proper} subformulas $\psi$ and $\omega$ (that is, all subformulas except for those that
	are only contained in $\sigma$ as part of helper formulas), then it is $\PEL1$-compatible. We will denote the set
	of helper formulas as $\Delta$.

	Let us assume that $\sigma$ is a theorem in $\ML$.

	By the subformula property, there is a derivation of $\sigma$ in $\ML$ that uses only the subformulas of $\sigma$.

	It is possible to translate this derivation to obtain a derivation of $t(\phi)$ in $\PEL1$.

	All the $\ML$ inference rules except for $\to\texttt{I}$ are also $\PL$ (and thus $\PEL1$) rules, so the corresponding steps
	are translated without changes. By the $\texttt{Premise inflation}$, we could add all $\Delta$ to the antecedents of
	these steps.

	Now, consider any of the remaining steps of the form $\infer{\Gamma \vdash \eb{\psi \to \omega}}{\Gamma, \psi \vdash \omega}$. We need to prove that, if
	$\Gamma, \Delta, \psi \vdash \omega$ is a theorem in $\PEL1$, then $\Gamma, \Delta \vdash \eb{\psi \to \omega}$ is also a theorem.

	Let us assume that both $\psi$ and $\omega$ are proper subformulas. From $\Gamma, \Delta, \psi \vdash \omega$ we obtain that
	$\psi \equ_{\Gamma,\Delta} \eb{\psi \land \omega}$. Now, once we substitute $\psi$ as $\phi$, $\psi \land \omega$ as $\psi$, $\psi \land \omega$ as $\chi$ and $\Gamma, \Delta$ as $\Gamma$ into the $\texttt{E}1$
	inference rule, we obtain $\eb{\psi \to \eb{\psi \land \omega}} \equ_{\Gamma,\Delta} \eb{\eb{\psi \land \omega} \to \eb{\psi \land \omega}}$. As $\eb{\eb{\psi \land \omega} \to \eb{\psi \land \omega}}$ is a helper
	formula, we obtain $\Gamma, \Delta \vdash \eb{\psi \to \eb{\psi \land \omega}}$. As $\eb{\eb{\psi \to \psi \land \omega} \to \eb{\psi \to \omega}} \in \Delta$, by $\to \texttt{E}$ rule we obtain
	$\Gamma, \Delta \vdash \eb{\psi \to \omega}$.

	If $\psi$ and $\omega$ are not proper subformulas, then, as $\psi \to \omega$ is a subformula of $\sigma$, it is either a helper
	formula (in which case $\Gamma, \Delta \vdash \eb{\psi \to \omega}$ is immediately obtained by $\texttt{x2x}$ and $\texttt{Premise inflation}$ rules),
	or a formula of a form $\psi' \to \psi' \land \omega$, where both $\psi'$ and $\omega$ are proper subformulas. In that latter case, by
	applying the reasoning for proper formulas above, we obtain $\Gamma, \Delta \vdash \eb{\psi' \to \eb{\psi' \land \omega}}$ as an intermediate
	step.

	By repeating the steps over and over, we obtain $\sigma$ as a theorem in $\PEL1$.

	Now, let us consider a mapping which turns sequent $\sigma = \Gamma \vdash \phi$ to the sequent $\Gamma, \Delta$, where $\Delta$ is
	the list of all helper formulas for $\sigma$. Such a mapping obviously turns theorems into theorems in every
	logic with $\texttt{Premise inflation}$ rule. It turns non-theorems into non-theorems in $\ML$, since it only adds
	as antecedents the consequent parts of antecedent-less theorems. As was shown above, it maps $\ML$
	sequents to $\PEL1$-compatible ones.

	Additionally, it takes polynomial time to obtain the resulting sequent.
}

\lemma{mlReducibleToPel2}
{
	$\ML$ is reducible to $\PEL2$ with respect to a polynomial-time mapping.
}
{
	The proof is similar to the previous one, except that we use

	$$\eb{\psi \to \psi}, \quad \eb{\eb{\psi \to \eb{\psi \land \omega}} \to \eb{\psi \to \omega}}$$

	helper formulas instead.

	By substituting $\psi$ as $\phi$, $\psi \land \omega$ as $\psi$, $\psi$ as $\chi$ and $\Gamma$ as $\Gamma$ into the $\texttt{E}2$ inference rule, we obtain
	$\Gamma, \Delta \vdash (\psi \to (\psi \land \omega))$ from $\Gamma, \Delta, \psi \vdash \omega$.

	Again, by the assumptions $\psi \to \omega$ and $\eb{\psi \to \eb{\psi \land \omega}} \to \eb{\psi \to \omega}$ we obtain $\Gamma, \Delta \vdash \psi \to \omega$.
}

\lemma{pelComplexity}
{
	$\ML$ and $\CL$ are reducible to $\PEL1$, $\PEL2$ and $\PEL$ with respect to a polynomial-time mapping.
	All three $\PEL$ logics are PSPACE-hard.
}
{
	$\PEL1$ and $\PEL2$ are subsystems of $\PEL$, which is, in turn, a subsystem of $\ML$. By Lemma \ref{middleLemma}, $\ML$
	is reducible to $\PEL$ by any of the mappings obtained above. As $\ML$ is PSPACE-hard, all the logics to
	which it is reducible with respect to the polynomial-time mapping should also be PSPACE-hard.

	$\CL$ is reducible to $\ML$, which is in turn reducible to each of the three $\PEL$s, and the mappings used
	are polynomial-time ones. Composing these mappings we obtain that $\CL$ is reducible to any of the
	$\PEL$s.
}

\subsection{PEL extensions complexity}

\theorem{pelExtensionsComplexity}
{
	Let $X$ be any logic containing $\PEL1$ (or $\PEL2$) and contained in $\CL$. The derivability
	problem for $X$ is co-NP-hard.
}
{
	It is enough to apply Lemma \ref{middleLemma} to the $\PEL1 < X < \CL$ (or $\PEL2 < X < \CL$) inequality. From
	the fact that $\CL$ is reducible to $X$ with respect to a polynomial-time mapping we immediately obtain
	that the derivability problem for $X$ is co-NP-hard.
}

\section{PEL\textsubscript{0} complexity}

\subsection{Framework}

For the sake of convenience, in this chapter we will redefine $\PEL_0$ as the logic obtained from $\PL$ by
adding the following inference rule:

\begin{longtable}{c c}
$\texttt{E}_0$ & $\infer{\eb{\phi_1 \to \psi_1} \vdash \eb{\phi_2 \to \psi_2}}{\phi_1 \vdash \phi_2 \qquad \phi_2 \vdash \psi_2 \qquad \psi_1 \vdash \psi_2 \qquad \psi_2 \vdash \psi_1}$
\end{longtable}

It is easy to see that this new inference rule is equivalent to the combination of $\texttt{E1}_0$ and $\texttt{E2}_0$.

We will define the set of \emph{significant} implications for the certain derivation of theorem $\sigma$ in $\PEL_0$
as the minimal set that:

\begin{itemize}
\item{contains all implication subformulas of $\sigma$}
\item{for each application of the $\texttt{E}_0$ rule, either contains both its left-hand and right-hand sides simultaneously or does not contain neither; and}
\item{contains all implication subformulas of all its elements.}
\end{itemize}

We will say that the formula $\omega$ is \emph{substantial} for the certain derivation step, if its structure matters
for the corresponding inference rule. We will say that all other formulas are \emph{auxiliary} for that step.
For example, if $\psi \vdash \phi \to \psi$ was obtained by applying $\to \texttt{IW}$ to $\psi \vdash \psi$ (which, in turn, was obtained by
$\texttt{x2x}$), both $\phi$ and $\psi$ are auxiliary formulas for this step, while $\phi \to \psi$ is significant.

\subsection{Elimination of insignificant formulas}

\lemma{existenceOfDerivationWithoutInsignificantImplications}
{
	For every $\PEL_0$ theorem $\sigma$ there is a derivation such that every its implication subformula
	is significant.
}
{
	Let us take any derivation of $\sigma$.

	Let $F$ be the set of all significant formulas for the derivation.

	Let us define the mapping $t_\sigma = t_F$ as follows:

	\begin{itemize}
	\item{$t_F\eb{x} = x$, where $x$ is a variable or constant}
	\item{$t_F\eb{\phi \land \psi} = \eb{t_F\eb{\phi} \land t_F\eb{\psi}}$}
	\item{$t_F\eb{\phi \to \psi} = \eb{t_F\eb{\phi} \to t_F\eb{\psi}} = \eb{\phi \to \psi}$ if $\eb{\phi \to \psi} \in F$, and $t_F\eb{\psi}$ otherwise}
	\item{$t_F\eb{\phi_1, \ldots, \phi_n \vdash \psi} = t_F\eb{\phi_1}, \ldots, t_F\eb{\phi_n} \vdash t_F\eb{\psi}$}
	\end{itemize}

	It is obvious that, for every $\phi \in F$, $t_F\eb{\phi} = \phi$. In particular, $t_\sigma\eb{\sigma} = \sigma$.

	Now, let us apply $t_F$ to all theorems of the derivation. We will prove that this application produced the derivation of $\sigma$ such that every its implication subformula is
	significant.

	It is obvious that every implication subformula of the result is significant for the original derivation. Additionally, every implication subformula of the result is significant
	for the result because for every $\phi \in F$, $\phi$ is left intact by $t_F$ and thus is the part of the result. Now it only remains to prove that the result is the correct derivation.

	Let us suppose that $\tau$ is the first incorrect derivation step of the result. Or, in other words, that $\tau$
	is the first sequent in the resulting pseudo-derivation such, that it could not be obtaining by applying
	some $\PEL_0$ inference rule to the preceding steps (which are all theorems).

	\begin{itemize}
	\item{
		It obviously could not originally be obtained by one of the rules $\top, \texttt{x2x}, \land\texttt{E}, \land\texttt{I}, \texttt{Cut}$, since every
		implication subformula is auxiliary for these rules
	}
	\item{
		It could not originally be obtained by the rule $\to\texttt{E}$ or $\to\texttt{IW}$, where $\eb{\phi \to \psi} \in F$, since both $\phi$
		and $\psi$ remain intact under $t_F$ transformation, and thus $\tau$ is obtained as in the original derivation
	}
	\item{
		It could not be originally obtained by the rule $\to\texttt{E}$, where $\eb{\phi \to \psi} \notin F$, since $t_F\eb{\Gamma \vdash \eb{\phi \to \psi}} =
		t_F\eb{\Gamma} \vdash t_F\eb{\psi}$; the transformation turns both the latter of the premises and the conclusion into
		the same sequent
	}
	\item{
		It could not be originally obtained by the rule $\to\texttt{IW}$, where $\eb{\phi \to \psi} \notin F$ by the same argument;
		the transformation turns both the premise and the conclusion into the same sequent
	}
	\item{
		It could not be originally obtained by the rule $\texttt{E}_0$, where $\eb{\phi_1 \to \psi_1} \in F$, since that would imply
		$\eb{\phi_2 \to \psi_2} \in F$ by the construction of $F$, and thus all of $\phi_1, \phi_2, \psi_1, \psi_2$ remain intact under $t_F$
		transformation, and thus $\tau$ is obtained as in the original derivation
	}
	\item{
		It could not be originally obtained by the rule $\texttt{E}_0$, where $\eb{\phi_1 \to \psi_1} \notin F$, since that would imply
		$\eb{\phi_2 \to \psi_2} \notin F$ by the construction of $F$, and thus $t_F\eb{\eb{\phi_1 \to \psi_1} \vdash \eb{\phi_2 \to \psi_2}} = t_F\eb{\psi_1 \vdash \psi_2}$;
		the transformation turns both the latter of the premises and the conclusion into the same sequent,
		yet the result does not contain duplicate sequents
	}
	\end{itemize}

	We just have exhausted all the possibilities for how the preimage of $\tau$ could be originally obtained.
	Thus our assumption of $\tau$ being the first incorrect derivation step is wrong; the result is indeed the correct
	derivation.
}

\subsection{Elimination of non-subformula implications}

\lemma{existenceOfDerivationWithoutNonsubformulaImplications}
{
	For every $\PEL_0$ theorem $\sigma$ there is a derivation such that every its implication subformula
	is a subformula of $\sigma$.
}
{
	Let us suppose that there is a theorem $\sigma$ such that every its derivation contains an implication
	subformula which is not a subformula of $\sigma$.

	For every derivation of $\sigma$, let us consider the pair $\eb{l, n}$, where $l$ is the length of the longest
	implication subformula of the derivation which is not a subformula of $\sigma$, and $n$ is the number of
	different implication $l$-length subformulas of the derivation which are not subformulas of $\sigma$.

	Note that every significant implication subformula of the derivation, which is not a subformula of
	$\sigma$, is an antecedent or consequent of a derivation step which was obtained by applying $\texttt{E}_0$ rule, or a
	subformula of such antecedent or consequent. Every significant implication formula of the derivation
	of the length $l$, which is not a subformula of $\sigma$, is an antecedent or consequent of a derivation step
	which was obtained by applying $\texttt{E}_0$ rule.

	Now among all derivations of $\sigma$ such that every their implication formula is significant, let us
	consider one with the smallest pair $\eb{l, n}$ (in a lexical order).

	Let us take the first step in this derivation such that it is obtained by the application of $\texttt{E}_0$ rule
	(thus its result has a form of $\eb{\phi_1 \to \psi_1} \vdash \eb{\phi_2 \to \psi_2}$), and either $\eb{\phi_1 \to \psi_1}$ or $\eb{\phi_2 \to \psi_2}$ is of the
	length $l$ and is not a subformula of $\sigma$. For the sake of simplicity let us suppose that $\phi_1 \to \psi_1$ is of the
	length $l$ and is not a subformula of $\sigma$; and that $\phi_1$ is longer than $\phi_2$.

	Note that at this step it is already established that $\phi_1 \vdash \phi_2$, $\phi_2 \vdash \phi_1$, $\psi_1 \vdash \psi_2$, and $\psi_2 \vdash \psi_1$ are
	theorems; and their derivations do not employ $\texttt{E}_0$ rules that would result in obtaining an implication
	of length $l$ which is not a subformula of $\sigma$.

	By transforming these derivations as in lemma \ref{existenceOfDerivationWithoutInsignificantImplications} we obtain derivations of $\phi_1 \vdash \phi_2$, $\phi_2 \vdash \phi_1$,
	$\psi_1 \vdash \psi_2$, and $\psi_2 \vdash \psi_1$ such, that every implication subformula of every (not only $\texttt{E}_0$) step of these
	derivations is either of the length less than $l$, or a subformula of $\sigma$. In particular, these derivation
	do not mention $\eb{\phi_1 \to \psi_1}$. Prepending these new derivations to the chosen derivation of $\sigma$ does not
	increase its $\eb{l, n}$ pair.

	Additionally, if in the transformed derivation there are any steps of the form $\phi_1 \vdash \phi_3$, $\phi_3 \vdash \phi_1$,
	$\Gamma \vdash \phi_1$, we append $\phi_2 \vdash \phi_3$, $\phi_3 \vdash \phi_2$, $\Gamma \vdash \phi_2$ respectively immediately after these (or after the
	corresponding equivalence was obtained), employing $\texttt{Cut}$ rule. As we didn’t introduce new applications
	there, such a modification does not change $\eb{l, n}$ pair of the derivation.

	Now let us replace every occurrence of $\eb{\phi_1 \to \psi_1}$ in the derivation with $\phi_2 \to \psi_1$. Let then us
	prove that, after replacement, derivation is still legitimate.

	Let us suppose $\tau$ is the first incorrect derivation step of the result.

	\begin{itemize}
	\item{
		It obviously could not be originally obtained by one of the rules $\top$, $\texttt{x2x}$, $\land\texttt{E}$, $\land\texttt{I}$ or $\texttt{Cut}$, since every
		implication subformula is auxiliary for these rules.
	}
	\item{
		It could not be originally obtained by $\to\texttt{E}$ or $\to\texttt{IW}$ rule, if $\eb{\phi_1 \to \psi_1}$ is auxiliary for these rules.
	}
	\item{
		It could not be originally obtained by $\to\texttt{E}$ rule, if $\eb{\phi_1 \to \psi_1}$ is a substantial for this rule,
		as we already have $\Gamma \vdash \phi_2$ at this moment.
	}
	\item{
		It could not be originally obtained by $\to\texttt{IW}$ rule, as $\Gamma \vdash \eb{\phi_2 \to \psi_1}$ is immediately obtained from
		$\Gamma \vdash \psi_1$ by the same rule.
	}
	\item{
		If could not be originally obtained by $\texttt{E}_0$ rule. The original step in that case should be of the
		form $\eb{\phi_1 \to \psi_1} \vdash \eb{\phi_3 \to \psi_3}$ or $\eb{\phi_3 \to \psi_3} \vdash \eb{\phi_1 \to \psi_1}$; yet both $\eb{\phi_2 \to \psi_1} \vdash \eb{\phi_3 \to \psi_3}$ and
		$\eb{\phi_3 \to \psi_3} \vdash \eb{\phi_2 \to \psi_1}$ could be obtained immediately by the same rule, as we already have
		$\phi_2 \equ \phi_3$.
	}
	\end{itemize}

	If $\phi_1$ was shorter than $\phi_2$, then $\psi_1$ must be longer than $\psi_2$. We perform the same operation
	then, except that we replace $\eb{\phi_1 \to \psi_1}$ with $\eb{\phi_1 \to \psi_2}$. The proof of the correctness of the changed
	derivation remain the same, except that, in $\to\texttt{IW}$ case, we obtain $\Gamma \vdash \eb{\phi_1 \to \psi_2}$ from $\Gamma \vdash \psi_2$.

	We just have constructed the new derivation of $\sigma$ having one implication subformula (which is not
	a subformula of $\sigma$) of the length $l$ less; thus either $l$ has decreased, or $l$ has remained the same while $n$
	has decreased. Yet we specifically considered the derivation with the smallest pair $\eb{l, n}$. This means
	that there is no such theorem $\sigma$; and that for every theorem there is a derivation containing only
	subformulas of $\sigma$ as implications.
}

\corollary{corollary71}
{
	If formulas $\phi_1, \ldots, \phi_n, \psi$ do not contain different equivalent proper subformulas, then
	the sequent $\sigma = \phi_1, \ldots, \phi_n \vdash \psi$ of $\PEL_0$ is $\PL$-compatible.
}
{
	Let $\sigma$ be a theorem of $\PEL_0$. Let us consider its derivation that does not contain implication
	subformulas that are not subformulas of $\sigma$.

	This derivation does not employ $\texttt{E}_0$ rule, as both parts of its conclusion are subformulas of $\sigma$, and
	therefore all parts of its premises are proper subformulas of $\sigma$, and therefore equivalent parts of its
	premises are actually equal, and therefore both parts of its conclusion are equal, and could as well be
	obtained using $\texttt{x2x}$ rule. Therefore, this derivation is a derivation in $\PL$.
}

\subsection{Elimination of equivalent subformulas}

We say that the set of formulas $\psi_i$ is \emph{free of equivalents} if they do not contain different subformulas
which are equivalent in $\PEL_0$.

\lemma{lemma73}
{
	Given set of formulas $\phi_i$ with the combined length of $n$, it is possible to compute, in
	$O\eb{n^3}$ time, the set of formulas $\psi_i$ free of equivalents such that, for every $i$, $\phi_i \equ \psi_i$ in $\PEL_0$; and $\psi_i$
	is no longer than $\phi_i$.
}
{
	Let us consider expression trees corresponding to the formulas $\phi_i$.

	Now, let us mark and transform these trees in several steps as follows:

	In the beginning, all the nodes are unmarked.

	On every step, we will consider the shortest (in terms of the length of the corresponding formula)
	unmarked node. If there is any marked node equivalent to the one under consideration, we’ll replace
	the latter with the former. In any case, we’ll mark the resulting node.

	We will continue doing such steps until all nodes are marked.

	It is easy to see that the following invariants are satisfied on every step:

	\begin{itemize}
	\item{
		Every ancestor (subformula) of the marked node is marked. This is obviously true in the beginning (since
		there are no marked nodes), and remains true on every step (since we only mark either shortest
		unmarked node, or the one which already satisfies this condition).
	}
	\item{
		Set of formulas represented by the marked nodes is free of equivalents. This is obviously true
		in the beginning, and remains true on every step (since we only mark a new node if formula it
		represents is not equivalent to any of formulas represented by the marked nodes, and all of its
		ancestors are already marked).
	}
	\item{
		Any unmarked node is not shorter than any marked node (in terms of the length of the corresponding subformulas). This is obviously true in the beginning, and remains true on every step
		(since we only mark a new node if it is shortest among the unmarked nodes).
	}
	\item{
		For every tree, the formula it represents is equivalent to, and not longer than the formula it
		originally represented (since every change in its structure is caused by replacing unmarked node
		with one not longer and equivalent to it).
	}
	\end{itemize}

	Thus, in the end we get the set of formulas $\psi_i$ satisfying the condition of lemma.

	It remains to estimate computational complexity of the process.

	There are at most $n$ steps, since every step reduces the number of unmarked nodes by one, and
	there was at most $n$ unmarked nodes.

	On every step, we have to find the shortest node ($O\eb{n}$), and then, for every marked node (number
	of which is less than $n$), to check whether the corresponding formulas are equivalent in $\PEL_0$. As all
	their ascendants are already marked (and thus corresponding subformulas are free of equivalents), by
	corollary \ref{corollary71} the check could be performed in $\PL$ in linear time. Thus, the step could be performed in
	$O\eb{n^2}$, and all the process could be performed in $O\eb{n^3}$ time.
}

\subsection{Multiple derivability problem complexity}

\lemma{multipleDerivabilityPel0}
{
	Multiple derivability problem for $\PEL_0$ is solvable in cubic time.
}
{
	Given hypotheses and queries, we could obtain the equivalent hypotheses and queries, free
	of equivalents, such that the combined length of new hypotheses and queries is not greater than the
	combined length of the original hypotheses and queries; the construction of new hypotheses and queries
	would require cubic time, by lemma \ref{lemma73}. Thus, multiple derivability problem for original hypotheses
	and queries in $\PEL_0$ is equivalent to the multiple derivability problem for new hypotheses and queries
	in $\PL$, which is, in turn, is solvable in linear time.
}

\corollary{corollary72}
{
	$\PEL_0 < \PEL1$, $\PEL_0 < \PEL2$
}
{
	As $\PEL1$ and $\PEL2$ are co-NP-hard, while $\PEL_0$ is cubic-time, $\PEL1$ and $\PEL2$ are not equal to
	$\PEL_0$.
}

\section{Logics with disjunction}

All of the results above could be directly translated to the case of logics with disjunction, with the
disjunction defined as follows:

\begin{longtable}{c c}
$\lor \texttt{E}$ & $\infer{\Gamma \vdash \chi}{\Gamma, \phi \vdash \chi &\quad \Gamma, \psi \vdash \chi &\quad \Gamma \vdash \eb{\phi \lor \psi}}$ \\\\
$\lor \texttt{I}_l, \lor \texttt{I}_r$ & $\infer{\Gamma \vdash \eb{\phi \lor \psi}}{\Gamma \vdash \phi},\quad \infer{\Gamma \vdash \eb{\phi \lor \psi}}{\Gamma \vdash \psi}$
\end{longtable}

The only diffence is that, as $\PL\withor$ is co-NP-hard, $\PEL_0\withor$ is co-NP-hard as well by lemma \ref{middleLemma}. As
a consequence, it is unclear whether $\PEL_0\withor$ is different from all of $\PEL1\withor$, $\PEL2\withor$, $\PEL\withor$, or not.

\section{Conclusion}

In this work we have defined the principle of equivalent formula substitution and the closure of primal
logic for this principle, and established the following facts:

\begin{enumerate}
	\item{$\CL\withor$ is reducible to $\CL$; $\CL$ is reducible to $\CL\withor$; $\CL$ is co-NP-complete.}
	\item{$\IL$ and $\ML$ are sound and complete with respect to the Kripke frames.}
	\item{$\IL$ is reducible to $\ML$. The same is true for their disjunction-free counterparts.}
	\item{$\ML$ is reducible to $\PEL1$ and $\PEL2$, and so the latter are PSPACE-complete. The same is true for
their disjunction-free counterparts.}
	\item{$\CL$ is reducible to any $\PEL1$ and $\PEL2$, which are therefore co-NP-hard.}
	\item{$\PEL_0$ is reducible to $\PL$; the derivability problem in $\PEL_0$ is solvable in cubic time}
\end{enumerate}

\end{document}